\documentstyle[a4]{article}

\input amssymb.sty
\input epsf
\def\si{\sigma}
\def\om{\omega}
\def\ip{\lrcorner}

\def\map{\rightarrow}
\def\bq{\begin{equation}}
\def\eq{\end{equation}}

\def\fg{{\frak g}}

\def\fp{{\frak p}}

\begin{document}
\section*{A note on generating functions}
Suppose $A$ is an affine symplectic space. There is an affine-invariant view of generating functions of symplectic transformations of 
$A$. Namely, let $H$ be a function on $A$. At any point $x$ we take the vector $u_x$ defined by $d_xH=u_x\ip \om$ ($\om$ is the 
symplectic form) and put it in $A$ so that $x$ lies in its middle. Then the map $\Phi_H$ sending the tails of $u_x$'s to their heads is a 
symplectic transformation:
$$\epsfxsize 8cm \epsfbox{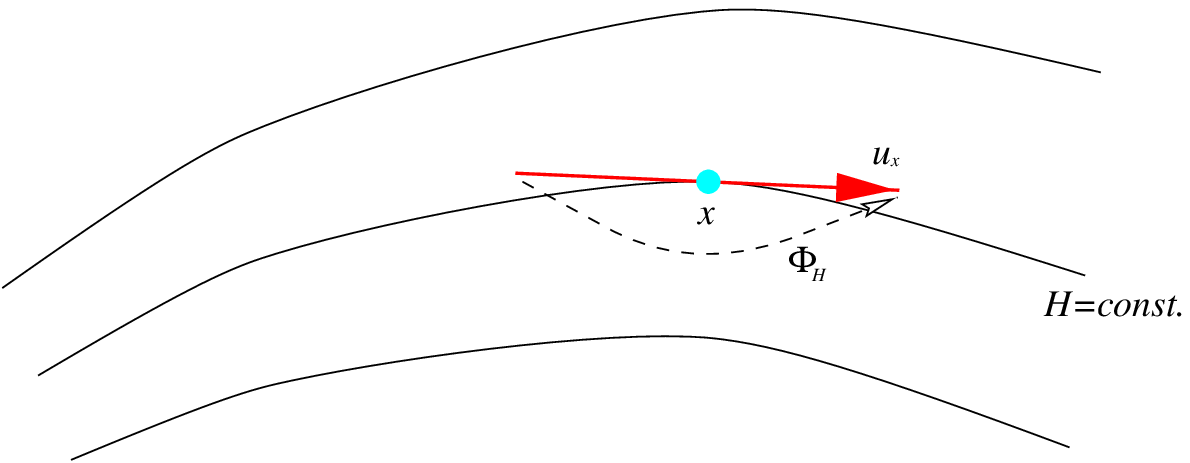} $$
Notice that for infinitesimal $H$, this is the usual infinitesimal transformation generated by Hamiltonian $H$. The map $H\mapsto\Phi_H$ 
is a kind of Cayley transform: choosing an origin in $A$ (to turn it to a vector space) and restricting ourselves to quadratic forms, we 
get the usual Cayley transform ${\frak {sp}}\map Sp$.

Symplectic transformations can be composed. The corresponding composition of generating functions is 
$H(x)=H_1(x_1)+H_2(x_2)+\mbox{symplectic area of }\triangle PQR $:
$$\epsfxsize 6cm \epsfbox{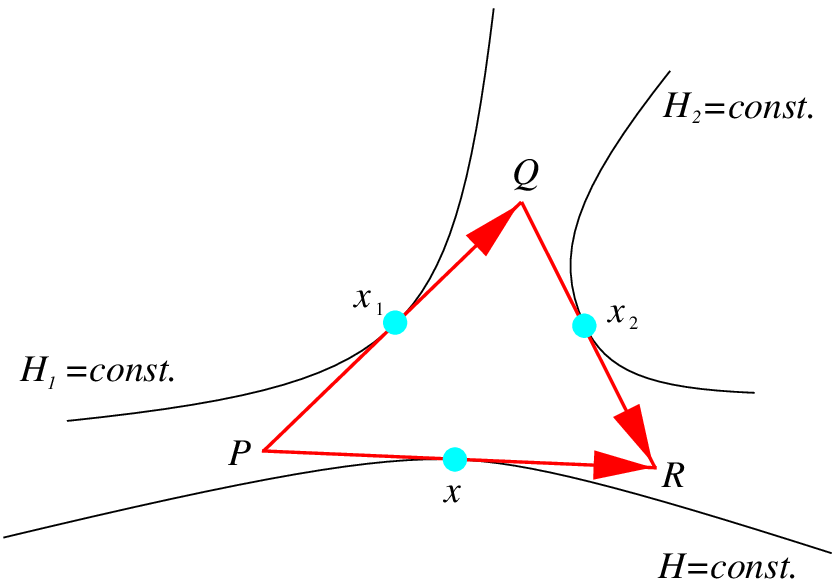}$$
Recall that the integral kernel of the Moyal product is $K(x_1,x_2,x)=\exp(\sqrt{-1}\times \mbox{symplectic area of }\triangle 
PQR/\hbar)$. We may notice that $\exp(\sqrt{-1}H/\hbar)$ is the classical part of the Moyal product of $\exp(\sqrt{-1}H_1/\hbar)$ and 
$\exp(\sqrt{-1}H_2/\hbar)$.

Let us have a look where these claims come from. A symplectic transformation of $A$ is (more-or-less) the same as a Lagrangian 
submanifold of $\bar A\times A$ (the graph of the map). For each point $x\in A$ the symmetry with respect to $x$ is a symplectic map. 
Identity is also a symplectic map, so that we have many Lagrangian submanifolds of $\bar A\times A$: $$\epsfxsize 4cm \epsfbox{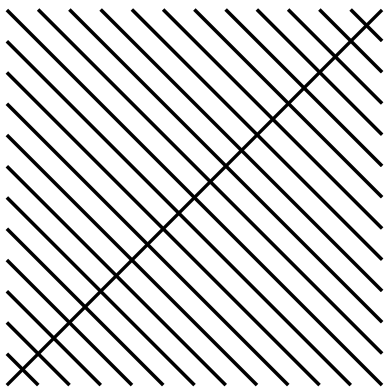}$$ 
In this way we have an isomorphism between $\bar A\times A$ and $T^*A$. Explicitely (as one immediately sees from the picture), a pair 
$(P,Q)\in\bar A\times A$ corresponds to $((P+Q)/2,(Q-P)\ip\om)\in T^*A$. Here the vector-and-its-midpoint picture appears.

Correspondence between generating functions and symplectic transformations is clear now: $dH$ is a Lagrangian submanifold of $T^*A$, and 
therefore of $\bar A\times A$. Let us also have a look where the composition law comes from. $\bar A\times A$ is a symplectic groupoid 
(the pair groupoid of $A$). The graph of its multiplication is a Lagrangian submanifold; using the identification of $T^*A$ and $\bar 
A\times A$, it should be given by a closed 1-form on $A\times A\times A$; this 1-form is the differential of the function 
$(x_1,x_2,x)\mapsto\mbox{symplectic area of }\triangle PQR $. The composition of generating functions and its connection with Moyal 
product follows.

For the fun of it, let us make a similar construction, replacing $A$ by the sphere $S^2$ with the area 2-form. Again, symmetry with 
respect to a point is a symplectic map, therefore we locally have a similar identification between $\overline{S^2}\times S^2$ and 
$T^*S^2$; more precisely, there is an isomorphism between the subset of covectors in $T^*S^2$ of length less than 2 and 
$\overline{S^2}\times S^2$ with erased pairs of antipodal points. Explicitely, to a non-antipodal pair $(P,Q)$ we associate a point in 
$TS^2$ (and thus, via $\om$, a point in $T^*S^2$) as on the picture:$$\epsfxsize 5cm \epsfbox{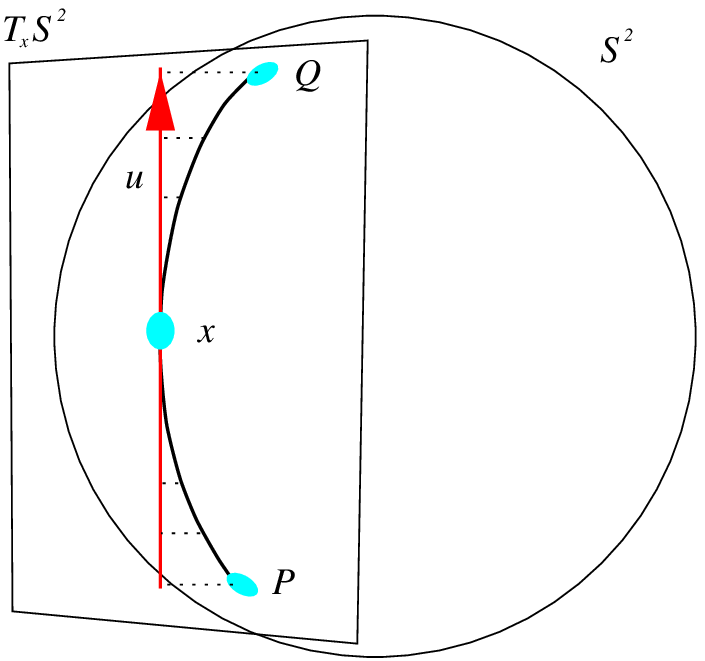} $$
$x$ is the midpoint of the shorter geodesic arc $PQ$ and $u\in T_xS^2$ appears by its orthogonal projection. This picture can be derived 
from the famous theorem of Archimedes, claiming that certain map between cylinder and sphere is area-preserving.

As a result, we have a similar picture of generating functions: for a function $H$ on $S^2$ and any point $x$ we take the vector $u_x$ 
defined by $d_xH=u_x\ip \om$, place it into the tangent plane so that $x$ is in its middle and project it into the sphere; $\Phi_H$ maps 
$P$ to $Q$. Composition rule looks as before (only triangles are spherical now).

Generally, this picture works with no changes for arbitrary symmetric symplectic space $M$. Using the symmetries we locally identify $\bar 
M\times M$ with $T^*M$. Multiplication in this pair groupoid is again given by the symplectic area of a surface bounded by the geodesic 
triangle $PQR$ with $x_1,x_2,x$ being the midpoints of its sides. The identification between $\bar M\times M$ and $T^*M$ is via a projection of $M$ into $T_xM$, as in the case $M=S^2$: Up to coverings, we embed $M$ into an affine space $A$. For any $x\in M$, the symmetry with respect to $x$ will be extended to an involution $\si_x$ of $A$; we project $M$ to $T_xM$  in the direction of $A^{\si_x}$ (the subspace of $A$ fixed by $\si_x$). Namely, since $M$ is a symmetric space, it is (a covering of) $G/G^\si$, where $G$ is a Lie group and $\si$ is an 
involutory automorphism of $G$. Let $\fg=\fg^\si\oplus\fp$ be the decomposition of $\fg$ to $\pm 1$ eigenspaces of $d\si$ (to make 
$G/G^\si$ into a symmetric symplectic space, one has to specify a $G^\si$-invariant symplectic form on $\fp$). As a homogeneous symplectic space, $M$ can be embedded (up to coverings) into an affine space over $\fg^*$ via (non-equivariant) moment map. If $x\in M$ is fixed by $G^\si$, $T_xM$ is  $\fp^*$ translated to $x$; we project $M$ 
to $T_xM$ in the direction of ${\fg^\si}^*$.

\vskip 3mm
{\it Pavol \v Severa, I.H.\'E.S. (IPDE postdoc)

severa@ihes.fr}

\end{document}